
\magnification=\magstep1
\input amstex
\UseAMSsymbols
\input pictex
\vsize=23truecm
\NoBlackBoxes
\pageno=1

   \font\rmk=cmr8    \font\itk=cmti8  

   \font\gross=cmbx10 scaled\magstep1 

\def\mo{\operatorname{mod}}

\def\Hom{\operatorname{Hom}}

\def\Ext{\operatorname{Ext}}

\def\bdim{\operatorname{\bold{dim}}}

\vglue1cm

\centerline{\gross Quiver Grassmannians for Wild Acyclic Quivers}
		          \bigskip

\centerline{Claus Michael Ringel}
		  	    \bigskip\bigskip
{\narrower\narrower \noindent Abstract.  A famous result of Zimmermann-Huisgen, Hille and Reineke asserts that any projective variety occurs as a 
quiver Grassmannian for a suitable representation of some wild acyclic quiver.
We show that this happens for {\it any} wild acyclic quiver. 
\par}

\plainfootnote{}{\rmk 2010 \itk Mathematics Subject Classification. \rmk
 16G20, 
 16G60, 
 14D20, 
 }

	\bigskip
Let $k$ be an algebraically closed field, and $Q$ a finite acyclic quiver.
The modules which we consider are the (finite-dimensional)
$kQ$-modules, where $kQ$ is the path
algebra of $Q$, thus the (finite-dimensional)
representations of $Q$ (with coefficients in $k$).  
We denote by $\mo kQ$ the corresponding module category. 

Let $M$ be a representation of $Q$ and $\bold d$ a dimension vector
for $Q$. The {\it quiver Grassmannian} $\Bbb G_{\bold d}(M)$ is the set 
of submodules of $M$ with dimension vector $\bold d$; this is a projective variety.
A famous result of Zimmermann-Huisgen, Hille and Reineke
asserts that any projective variety occurs as the quiver Grassmannian for some
wild acyclic quiver $Q$, see for example [R2]. 
In [R3] we have shown that one may take as quiver a
Kronecker quiver $Q = K(n)$, for $N$ a suitable reduced representation of $Q$ 
and as $\bold d$ the dimension vector $(1,1)$, see also [H]. Here, 
a representation of a Kronecker
quiver is called {\it reduced} in case it has no simple injective direct
summand. For a reduced representation $N$ of a Kronecker quiver, 
$\Bbb G_{(1,1)}(N)$ is the set of indecomposable submodules of $N$ 
of length 2. We call 
indecomposable modules of length 2 {\it bristles}, and $\beta(N) = 
\Bbb G_{(1,1)}(N)$ the {\it bristle variety} of $N$.
We use the result of [R3] in order to show:
   \medskip
{\bf Theorem.} {\it
If $Q$ be a wild acyclic quiver, then any projective variety occurs as a
quiver Grassmannian for a suitable representation of $Q$.}
    \medskip
For the proof, we will construct full exact subcategories $\Cal E$ of $\mo kQ$
which are equivalent to $\mo kK(n)$ with $n$ arbitrarily large. In order to
define such an $\Cal E$, we start with a pair $X,Y$ of orthogonal bricks with 
$\dim_k \Ext^1(Y,X) = n$, and $\Cal E = \Cal E(Y,X)$ will be the full
subcategory of all $kQ$-modules $M$ with an exact sequence of the
form 
$$
  0 @>>> X^a @>>> M @>>> Y^b @>>> 0,
$$
where $a,b$ are natural numbers. Always, $\bold x$ and $\bold y$ will denote
the dimension vectors of $X$ and $Y$, respectively. 
An equivalence between $\mo kK(n)$ and $\Cal E$ is given by
an exact fully faithful functor 
$$
  \eta\:\mo kK(n) \to \mo kQ
$$ 
with image $\Cal E$. 
We say that a module $M$ in $\Cal E$ is {\it $\Cal E$-reduced}
provided it has no direct summand isomorphic to $Y$, thus provided it is
the image of a reduced $kK(n)$-module under $\eta$. The 
indecomposable $kQ$-modules $U$ with an exact sequence of the form
$0 \to X \to U \to Y \to 0$ will be called $\Cal E$-bristles (of course, they
are the images under $\eta$ of the bristles in $\mo kK(n)$, note 
that $\Cal E$-bristles have dimension vector $\bold x+\bold y$).

For any $kK(n)$-module $N$, the functor $\eta$ identifies
the bristle variety $\beta(N)$ of $N$ 
with the set of submodules of $\eta N$ which are $\Cal E$-bristles. 
It remains to specify
conditions such that the set of $\Cal E$-bristles is just the quiver Grassmanian
$\Bbb G_{\bold x+\bold y}(\eta N)$. We will choose $X, Y$ so that the 
following closure condition (C) is satisfied:
	\medskip 
(C) {\it If $M$ is an $\Cal E$-reduced module in $\Cal E(Y,X)$
and $U$ is a submodule of $M$ with $\bdim U = \bdim X + \bdim Y$, then
$U$ is an $\Cal E$-bristle.}
	\medskip
If the condition (C) is satisfied, then for any reduced 
representation $N$ of $K(n)$, there is
a canonical bijection between $\Bbb G_{(1,1)}(N)$ and
$\Bbb G_{\bold x +\bold y }(\eta N)$.
Namely, if $B$ is a submodule of
the $kK(n)$-module $N$ with $\bdim B = (1,1)$, then $\eta B$ is a submodule of
$\eta N$ with dimension vector $\bold x+\bold y$. Conversely, 
if $U$ is a submodule
of $\eta N$ with $\bdim U = \bold x+\bold y$, then, by condition (C), $U$ belongs
to $\Cal E(Y,X)$, say $U = \eta B$ for some $K(n)$-submodule $B$ and the dimension vector of $B$ is $(1,1)$.
   \bigskip
Our aim is to exhibit for any wild acyclic quiver $Q$ 
and any natural number $m$
an orthogonal pair $X,Y$ of $kQ$-modules which are 
bricks such that $\dim_k\Ext^1(Y,X) = n \ge m$ and
such that the condition (C) is satisfied. The following well-known proposition 
suggests to deal with two different cases.
       \medskip 
{\bf Proposition.} {\it A wild acyclic quiver $Q$
with at least 3 vertices has a vertex $\omega$
which is a sink or a source such that the quiver $Q'$ obtained from $Q$ by deleting $\omega$
is connected and representation-infinite.} \hfill$\square$
	\bigskip 
{\bf Case 1.} Assume that $Q$ is a connected quiver with a vertex $\omega$
which is a sink or a source such that the quiver $Q'$ obtained from $Q$ by deleting $\omega$
is connected and representation-infinite.
Up to duality, we can assume that $\omega$ is
a source, thus there is an arrow $\omega \to p$ with $p\in Q'_0$.

Let $Y = S(\omega)$.
Since $Q'$ is connected and representation-infinite, there is an exceptional 
$kQ'$-module $X$
with $\dim_k X_p \ge m$. The arrow $\omega \to p$ shows that $\dim_k \Ext^1(Y,X) \ge \dim_k X_p$.
This pair $X,Y$ is the orthogonal pair of bricks which we use in order to look at
$\Cal E(Y,X)$.

	\medskip
{\bf Lemma 1.} {\it Let $a$ be a natural number. Any submodule $W$ of $X^a$ 
with $\bdim W = \bold x$ is isomorphic to $X$.}
	\medskip
Proof. We denote by $\langle-,-\rangle$ the bilinear form on the
Grothendieck group $K_0(kQ)$ with $\langle\bdim M,\bdim M') =
\dim_k \Hom(M,M') -\dim_k \Ext^1(M,M').$
Since $X$ is exceptional, we have $\langle X,W\rangle =
\langle X,X\rangle > 0$,  Therefore, there is a non-zero 
homomorphisms $f\:X \to W$. Let $\iota\:W \to X^a$ be the inclusion map.
The composition $\iota f\:X \to X^a$ is nonzero.  
Since $X$ is a brick, we see
that $f\:X \to W$ is a split monomorphism, in particular injective. Now
$\bdim X = \bdim W$ implies that $f$ is an isomorphism. 
\hfill$\square$
	\medskip
Proof of condition (C). Let $M$ be an $\Cal E$-reduced $kQ$-module
in $\Cal E(Y,X)$, say with an exact sequence
$$
  0 @>>> X^a @>\mu>> M @>\pi>> Y^b @>>> 0.
  $$
  Let $U$ be a submodule of $M$ with dimension vector $\bold x+\bold y$ and inclusion
  map $\iota\:U \to M$. The composition $\pi\iota$ is non-zero, since otherwise $U$
  would be a submodule of $X^a,$ but $\dim_k U_\omega = 1$ whereas $X_\omega = 0$.
If follows that the image of
  $\pi\iota$ is isomorphic to $Y$. If we denote the kernel of $\pi\iota$ by $W$,
  we obtain the following commutative diagram with exact rows and vertical monomorphisms:
  $$
\CD
    0 @>>> W @>>> U @>>> Y @>>> 0   \cr
     @.     @VVV         @V\iota VV      @VVV        \cr
       0 @>>> X^a @>\mu>> M @>\pi>> Y^b @>>> 0.
\endCD
$$
Of course, $\bdim W = \bold x$, thus Lemma 1 shows that $W$ is isomorphic to
$X$. In particular, $U$ belongs to $\Cal E$. 

It remains to show that $U$ is indecomposable. Otherwise,
$U$ would be isomorphic to $W\oplus U$. Thus $M$ would have a submodule 
isomorphic to $Y$. But $Y$ is relative injective inside $\Cal E$, 
thus $M$ would have a direct
summand isomorphic to $Y$, in contrast to our assumption that $M$
is $\Cal E$-reduced. This shows that $U$ is indecomposable, thus 
an $\Cal E$-bristle.
\hfill$\square$
    \bigskip
{\bf Case 2.} Here we consider the $3$-Kronecker quiver $Q = K(3)$, with
two vertices $1$ and $2$ and three arrows $\alpha,\beta,\gamma:1\to 2.$
Let $\lambda_1,\dots,\lambda_n$ be pairwise different non-zero elements of $k$ with $n\ge 2$.
Let $X = X(\lambda_1,\dots,\lambda_n) = (k^n,k^n;\alpha,\beta,\gamma)$ be defined by
$$
  \alpha(e(i)) = e(i),\quad \beta(e(i)) = \lambda_i e(i),\quad \gamma(e(i)) 
  = e(i+1),
$$
for $1\le i \le n$, 
where $e(1),\dots, e(n)$ is the canonical basis of $k^n$ and $e(n+1) = e(1).$
Let $Y = (k,k;1,0,0)$.
We denote by $Q'$ the subquiver of $Q$ with arrows $\alpha,\beta$, this is
the $2$-Kronecker quiver $K(2)$. 
For the structure of the module category of the $2$-Kronecker quiver $K(2)$, see
for example [R1]. 
The restriction of $X, Y$ to $Q'$ shows that $\Hom(X,Y) = \Hom(Y,X) = 0.$
The endomorphism ring of $X|Q'$ is $k\times\cdots\times k$; and the only endomorphisms
of $X|Q'$ which commute with $\gamma$ are the scalar multiplications. This shows
that $X$ is a brick. Also, it is easy to see that $\dim_k\Ext^1(Y,X) = n.$ 

   \medskip
{\bf Lemma 2.} {\it Let $a$ be a natural number. Any submodule $W$ of $X^a$ 
with $\bdim W$ of the form $(w,w)$ is isomorphic to $X^s$ for some $s$.}
	\medskip
Proof: Let $M = X^a$ and decompose $M|Q' = \bigoplus_{i=1}^n M(i)$, where 
$\beta(x) = \lambda_i x$ for $x\in M(i)_1$. Here, we use $\alpha$ in order to identify $M_1$ and $M_2$.
Now we consider the submodule $W$ of $M$. Note that $W|Q'$ has to be regular, since it cannot
have any non-zero preinjective direct summand. As a regular submodule of a semisimple regular
Kronecker module it has to be a direct summand of $M|Q'$, thus we have a similar direct decomposition $W = \bigoplus W(i)$, where $W(i) = W\cap M(i)$.

The linear map $\gamma$ restricted to $W(i)_1$
is a monomorphism $W(i)_1 \to W(i+1)_2 = W(i+1)_1$ for $1\le i \le n$;
we obtain in this way a monomorphism $W(1)_1 \to W(1)_2 = W(1)_1.$ This shows that all the
monomorphisms $W(i)_1 \to W(i+1)_2 = W(i+1)_1$ are actually bijections. 
Let $\dim_k W(1)_1 = s.$ It follows that $W$ is isomorphic to $X^s$. 
\hfill$\square$
	\bigskip
Proof of condition (C). Let $M$ be an $\Cal E$-reduced $kQ$-module in $\Cal E$
and let $U$ be a submodule of $M$ with dimension vector $\bold x + \bold y = (n+1,n+1)$ and with inclusion map $\iota\:U \to M$.
	
Starting with the exact sequence $0 @>>> X^a @>\mu>> M @>\pi>> Y^b @>>> 0$
and the inclusion map $\iota\:U \to M$, let 
$W$ be the kernel and $\overline U$ the image of $\pi\iota\: 
U \to Y^b$. We obtain the following 
commutative diagram with exact rows and injective vertical maps:
$$
  \CD
    0 @>>> W @>>> U @>>> \overline U @>>> 0   \cr
     @.     @VVV         @V\iota VV      @VVV        \cr
    0 @>>> X^a @>\mu>> M @>\pi>> Y^b @>>> 0;
  \endCD
$$
Let us consider the restriction of these modules to $Q'$. 
Since $M|Q'$ is regular, it has no non-zero preinjective direct summand.
Thus any submodule of $M|Q'$ with dimension vector $(n+1,n+1)$ has to be
regular. This shows that $U|Q'$ is regular. Actually, $M|Q'$ is semisimple
regular, thus also its regular submodule $U|Q'$ is semisimple regular
(and a direct summand of $M|Q'$). Next, $\pi\iota$ is a map between regular
$kQ'$-modules, it follows that the kernel $W|Q'$ and the image $\overline U|Q'$ 
are regular $kQ'$-modules. In particular, the dimension vector of
$W$ is of the form 
$\bdim W = (w,w)$ for some $0\le w \le n+1$.

Now $U|Q'$ is a regular submodule of the semisimple regular $kQ'$-module
$Y^b|Q'$, thus $\overline U|Q'$ is a 
direct sum of copies of $Y|Q'$.  
By construction, $Y$ is annihilated by $\gamma$. Since
$\overline U$ is a submodule of $Y^b$, it follows that $\overline U$ is
annihilated by $\gamma$. Altogether, we see that $\overline U$ is the
direct sum of copies of $Y$. 

We claim that $W\neq 0$. Otherwise $U = \overline U = Y^{n+1}$, thus
$Y$ is a submodule of $M$. But $Y$ is relative injective in $\Cal E$,
thus $Y$ would be a direct summand of $M$. However, by assumption, $M$
is $\Cal E$-reduced. This contradiction shows that $W\neq 0.$

Now $W$ is a submodule of $X^a$ with dimension vector $(w,w)$, thus,
according to Lemma 2, $W$ is a direct summand of say $s$ 
copies of $X$ and $s\ge 1$. The equality $(w,w) = (sn,sn)$ implies that
that $s = 1$, since $w \le n+1$
and $n\ge 2$. Is this way, we see that $W$ is isomorphic to $X$.
It follows that $\bdim \overline U = (1,1)$ and therefore $\overline U = Y$.

Finally, as in Case 1, we see that $U$ is indecomposable, using again the
assumption that $M$ is $\Cal E$-reduced. This shows that $U$ is an
$\Cal E$-bristle.
\hfill$\square$
	\bigskip
{\bf Remark.} We should stress that given orthogonal bricks $X,Y$
in $\mo kQ$, the condition (C) is usually not satisfied. Here is a typical example
for $Q = K(3)$. As above, let $Y = (k,k;1,0,0)$, but for $X$ we now take 
$X = X'(\lambda_1,\lambda_2) = (k^2,k^2;\alpha,\beta,\gamma),$ defined by
$$
  \alpha(e(i)) = e(i),\quad \beta(e(i)) = \lambda_i e(i),\quad \gamma(e(1)) 
  = e(2),\quad \gamma(e(2)) = 0
$$
for $1\le i \le 2$. Again, $e(1),e(2)$ is the canonical basis of $k^2$ and
$\lambda_1 \neq \lambda_2$ are assumed to be non-zero elements of $k$.
Since $\dim_k\Ext^1(Y,X) = 2,$ there is an equivalence 
$\eta\:\mo kK(2)\to \Cal E(Y,X)$. Let $N$ be an indecomposable
$kK(2)$-module with dimension vector $(2,b)$ 
(note that $b$ has to be equal to $1,2$ or $3$) and $M = \eta N$. Thus
there is an exact sequence
$$
    0 @>>> X^2 @>>> M @>>> Y^b @>>> 0.
$$ 
Since we assume that $N$ is indecomposable, it is reduced, thus
$M$ is $\Cal E$-reduced. Note that $X$ has a (unique) 
$kQ$-submodule $V$ with dimension vector $(1,1)$: 
the vector spaces $V_1$ and $V_2$ both are generated by $e(2)$. 
The submodule $U = X\oplus V$ of $X^2$ is a submodule of $M$ with
dimension vector $(3,3) = \bold x + \bold y$, and it 
is not an $\Cal E$-bristle. Thus, condition (C) is not satisfied. Here, 
$\eta$ defines a proper embedding of $\beta(N) = \Bbb G_{(1,1)}(N)$ 
into $\Bbb G_{\bold x+\bold y}(M).$
	\bigskip
{\bf References}
     \medskip
\item{[H]} L\. Hille: Moduli of representations, quiver Grassmannians 
   and Hilbert schemes. arXiv: 1505.06008.

\item{[R1]} C\. M\. Ringel: Tame algebras and integral quadratic forms. 
   Springer Lecture Notes in Math\. 1099 (1984).
\item{[R2]} C\. M\. Ringel: 
    Quiver Grassmannians and Auslander varieties for wild algebras.
    J\. Algebra 402 (2014), 351-357. 
\item{[R3]} C\. M\. Ringel: The eigenvector variety of a matrix pencil.
  Linear Algebra and Appl. (to appear).  arXiv:1703.04097
	    \bigskip

{\rmk

Claus Michael Ringel \par
Fakult\"at f\"ur Mathematik, Universit\"at Bielefeld\par
D-33501 Bielefeld, Germany\par
E-mail: ringel\@math.uni-bielefeld.de
	\medskip

\bye